\documentclass[10pt,reqno]{article}
\usepackage{amsmath, amsthm, amssymb, stmaryrd}
\usepackage{hyperref}
\usepackage{enumerate}
\usepackage{url}
\usepackage{color}
\usepackage{tikz}
\usepackage{xcolor}
\usepackage{cancel}
\usetikzlibrary{patterns}
\usepackage[T1,T2A]{fontenc}
\usepackage[utf8]{inputenc}

\theoremstyle{definition}

\theoremstyle{remark}

\def \le {\leqslant}
\def \ge {\geqslant}

\usepackage{tikz} 
\usetikzlibrary{calc,arrows,decorations.pathreplacing,fadings,3d,positioning}

\begin{document}

\begin{Large}
\centerline {A note on  uniform version of  Littlewood inequality and Fibonacci numbers}
\end{Large}

\vskip+0.3cm
\begin{Large}
\centerline{by Nikolay Moshchevitin }
\end{Large}
\vskip+1cm

{\bf Abstract}: We give a simple proof of a recent result by Schleischitz \cite{S} dealing with a 
counterexample to  the uniform Littlewood conjecture. Our construction is based on simple properties of Fibonacci numbers.

\vskip+2cm
{\bf 1.  Brief introduction.}
\vskip+0.3cm

Answering a question from the paper \cite{K}, in \cite{S} the following result was proven:
{\it  there exist real numbers $\xi,\zeta$ such that 
$$
\limsup_{Q\to \infty} \, Q\, \left( \min_{x\in \mathbb{Z}_+:\,  1\le x < Q}
||x\xi|| \cdot ||x\zeta||\right) \ge 0.005326.
$$
where  $||\cdot ||$ denotes the distance to the nearest integer.
}

The proof from \cite{S} is based on a deep and fundamental result by Bourgain and Kontorovich \cite{B1,B2} and uses Vinogradov's inequality for  certain exponential sums.
Here we give an easy  proof of a more precise statement which is based on 
simple properties of Fibonacci numbers 
 only.

\vskip+0.3cm
{\bf 2. Uniform distribution lemma and special numbers.}
\vskip+0.3cm

Let $ F_n$ denotes the $n$-th Fibonacci number.
The following lemma used standard argument from uniform distribution theory and Möbius inversion.

\vskip+0.3cm
{\bf Lemma 1.} {\it
For any two intervals
$$
I = [u, u +\eta],\,\,\,\, J = [v, v+\eta] \subset [0,1]
$$
for all $n$ large enough, there exists an integer $a$ satisfying
$$
\frac{a}{F_n} \in I,\,\,\,\,\,
\left\{ \frac{F_{n-1}a}{F_n}\right\} \in J,\,\,\,\,\,1\le a < F_n,\,\,\,\,(a, F_n) = 1ю
$$
}
\vskip+0.3cm

We give a proof of Lemma 1 in Section 6 below.
Now we construct  real numbers $\alpha, \beta$ which admit special rational approximations.

\vskip+0.3cm

Let $ n_0 = 1$, $\delta_0 = 1$ and  $\delta_\nu \to 0$ be an arbitrary sequence of positive numbers.
Take
$$
I_0 = J_0 = [0,1]
$$
and apply the following inductive procedure. If  for an  integer $n_\nu$ 
numbers
$$
\alpha_\nu = \frac{a_\nu}{F_{n_\nu}},\,\,\,\,\,\,\
\beta_\nu =  \left\{\frac{F_{n_\nu -1}a_\nu}{F_{n_\nu}}\right\},\,\,\,\,\,
\text{where}\,\,\,\,\,
a_\nu\,\,\,\,\,\,
\text{is an integer},\,\,\,\,\,
(a_\nu, F_{n_\nu}) = 1
$$
are constructed, consider 
 segments
$$
I_\nu = \left[\alpha_\nu, \alpha_\nu +\frac{\delta_\nu}{F_{n_\nu}^2}\right],\,\,\,\,\,\,\,\,\,\,
J_\nu = \left[\beta_\nu, \beta_\nu +\frac{\delta_\nu}{F_{n_\nu}^2}\right].
$$
Take $n_{\nu+1}$ large enough and apply Lemma 1 to get  integer $a_{\nu+1}$ such that 
for rationals
$$
\alpha_{\nu+1} = \frac{a_{\nu+1}}{F_{n_{\nu+1}}},\,\,\,\,\,\,\
\beta_{\nu +1}= \left\{\frac{F_{n_{\nu+1} -1}a_{\nu+1}}{F_{n_{\nu+1}}}\right\}
$$
for the corresponding segments one has
$$
I_{\nu +1}\subset I_\nu ,\,\,\,\,\, J_{\nu +1}\subset J_\nu .
$$
Then we get numbers
$$
\alpha = \bigcap_\nu I_\nu,\,\,\,\,\, \beta = \bigcap_\nu J_\nu
$$
such that 
\begin{equation}\label{o1}
\max \left( 
\left|\alpha - \frac{a_\nu}{F_{n_\nu}}\right| ,
\left|\beta - \frac{F_{n_{\nu}-1}a_\nu}{F_{n_\nu}}\right|
\right)\le \frac{\delta_\nu}{F_{n_\nu}^2}\,\,\,\,\,
\text{for all}\,\,\,\,\, \nu.
\end{equation}

\vskip+0.3cm 

\vskip+0.3cm
{\bf 3. Result.}
\vskip+0.3cm

We prove the following
\vskip+0.3cm

{\bf Theorem 1.} {\it For real numbers $\alpha, \beta $  defined in the previous section one has
$$
\limsup_{Q\to \infty}\, Q\, \left( \min_{x\in \mathbb{Z}_+:\,  1\le x < Q}\, 
   ||\alpha x|| \cdot  ||\beta x||\right) \ge \frac{2}{3+\sqrt{5}} = 0.381966^+.
$$
}
\vskip+0.3cm

In the next section we will recall simple approximation properties of Fibonacci numbers and in Section 5 we finalise the proof of Theorem 1.

\vskip+0.3cm
{\bf 4. Approximations with Fibonacci numbers.}
\vskip+0.3cm
First of all we observe that if a rational fraction $\frac{y}{x}$ is not a convergent to  the number $\frac{F_{n-1}}{F_n}$, then
\begin{equation}\label{q1}
\left| \frac{F_{n-1}}{F_n} -\frac{y}{x}\right| \ge \frac{1}{2x^2} .
\end{equation}
Then,   if 
$\frac{y}{x} = \frac{F_{k-1}}{F_k}$, $ k <n$ is a convergent to
$$\frac{F_{n-1}}{F_n}
= [0;\underbrace{1,...,1}_{n-1}]б
$$
by the  approximation formula 
for any positive $\varepsilon$ and for $n $
large enough we have 
\begin{equation}\label{q2}
\left| \frac{F_{n-1}}{F_n} -\frac{F_{k-1}}{F_k}\right| =  
\frac{1}{F_k^2([\underbrace{1;1,...,1}_{n-k}]+ [0;\underbrace{1,...,1}_{k-1}]) }\ge 
\frac{1}{F_k^2\left(\frac{1+\sqrt{5}}{2}+\varepsilon+ 1\right) },
\end{equation}
as
$$
[\underbrace{1;1,...,1}_{n}] \to  \frac{1+\sqrt{5}}{2}\,\,\,\,\text{as}\,\,\,\,\,
n\to \infty,\,\,\,\,\,\,\,\,\,\,
\text{and}\,\,\,\,\,\,\,\,\,\,\
[0;\underbrace{1,...,1}_{n}]\le 1\,\,\,\,\,\, \text{for every}\,\,\,\,\,n.
$$
Combining (\ref{q1}) and (\ref{q2}) we see that 
$$
  \left|\left|
\frac{x}{F_n} \right|\right|  
\cdot
\left|\left|
\frac{F_{n-1}x}{F_n} \right|\right| =
\frac{|x|}{F_n} 
\cdot
\left|\left|
\frac{F_{n-1}x}{F_n} \right|\right| 
 \ge
\frac{1}{F_n\left(\frac{3+\sqrt{5}}{2}+\varepsilon\right) }\,\,\,\,\,
\text{for any}\,\,\,\,\, x\,\,\,\,\,
\text{in the range}\,\,\,\,\,\, 0<|x|\le \frac{F_n}{2},
$$
or
$$
 \min_{1\le x< F_n} \,\,
   \left|\left|
\frac{x}{F_n} \right|\right| 
\cdot
 \left|\left|
\frac{F_{n-1}x}{F_n} \right|\right| \ge
\frac{1}{F_n\left(\frac{3+\sqrt{5}}{2}+\varepsilon\right) }.
$$

Since,  for $a$ coprime to $F_n$ the numbers $ ax$ run through all the residues modulo $F_n$ when $x$ runs through all the residues modulo $F_n$ we conclude that 
\begin{equation}\label{q5}
 \min_{1\le x< F_n} \,\,
  \left|\left|
\frac{ax}{F_n} \right|\right| 
\cdot
  \left|\left|
\frac{F_{n-1}ax}{F_n}  \right|\right|\ge
\frac{1}{F_n\left(\frac{3+\sqrt{5}}{2}+\varepsilon\right) }\,\,\,\,\,
\text{for any}\,\,\,\,\, a\,\,\,\,\,\text{coprime to }\,\,\,\,\, F_n.
\end{equation}

\newpage
{\bf 5. Proof of Theorem 1.}
\vskip+0.3cm

For $\alpha,\beta$ constructed in Section 2  and for 
$x$ from the range $ 1\le x < F_{n_\nu}$ we consider the value
$$
||\alpha x||\cdot ||\beta x|| =
 \left|\left|
 \frac{a_{\nu}x}{F_{n_\nu} }+ \frac{\delta_{1,\nu}(x)}{F_{n_\nu}}
 \right|\right| 
 \cdot
 \left|\left|
 \frac{a_{\nu}F_{n_\nu-1}x}{F_{n_\nu} }+ \frac{\delta_{2,\nu}(x)}{F_{n_\nu}}
 \right|\right| =
  \left|\left|
 \frac{a_{\nu}x}{F_{n_\nu} }
 \right|\right| 
 \cdot
 \left|\left|
 \frac{a_{\nu}F_{n_\nu-1}x}{F_{n_\nu} }
 \right|\right| +\frac{\delta_{3,\nu}(x)}{F_{n_\nu}},
$$
where 
$$
\max_{1\le x <F_{n_\nu}}\,
 |\delta_{j,\nu} (x)|  \to 0,\,\,\,\,\, j = 1,2,3.
 $$
Now we apply (\ref{q5}) and 
continue with 
$$
||\alpha x||\cdot ||\beta x|| \ge
\frac{1}{F_{n_\nu}\left(\frac{3+\sqrt{5}}{2}+\varepsilon_1\right) }\,\,\,\,\,
\text{for all}\,\,\,\,\, x\,\,\,\,\,\text{in the range}\,\,\,\,\,
1\le x<F_{n_\nu}
$$
with small $\varepsilon_1$.
Everything is proven.$\Box$

\vskip+0.3cm
{\bf 6. Proof of Lemma 1.}
\vskip+0.3cm

 As fractions $\frac{F_{k-1}}{F_k}, k< n$ are convergent fractions to $\frac{F_{n-1}}{F_n}$  we see that
 $$
 \left|\left|\frac{F_{n-1}}{F_n} F_k\right|\right|  \le \frac{1}{F_k^2}.
 $$
 By a well known property (see, for example, Lemma 3 from Ch. I, \S 3 from \cite{Ca})   for any $\sigma>0$
 there exists $k_* = \frac{n}{2} +O(n^\sigma)$ such that $ (k_*,n) =1 $ and so
 $$
 (F_{k_*}, F_n) = F_{(k_*,n)} = 1.
$$ 
 As the sequence 
$$
\left\{
\frac{F_{n-1}}{F_n}x\right\}, \,\,\,\,\,\,\,\,\, uF_n \le x\le \left(u+\frac{\eta}{2}\right)F_n
$$
is uniformly distributed (moreover, it has logarithmic discrepancy, see Ch. 2 \S3 from \cite{KN}) we conclude that  there exists $a_0$
satisfying
$$
\left\{
\frac{F_{n-1}}{F_n}a_0\right\} \in
\left(
v+\frac{\eta}{3}, v+\frac{2\eta}{3}\right)
 ,\,\,\,\,\,\text{and}\,\,\,\,\,
uF_n \le a_0\le \left(u+\frac{\eta}{2}\right)F_n.
$$
Again,  by Lemma 3 from Ch. I, \S 3 from \cite{Ca}, for any $\sigma>0$  there exists $j = O_\sigma(F_n^\sigma)$ such that 
$$
a = a_0 + j F_{k_*}\,\,\,\,\, \text{is coprime to}\,\,\,\,\, F_n.
$$
Now
$$
\left\{
\frac{F_{n-1}}{F_n}a\right\}=
\left\{
\frac{F_{n-1}}{F_n}(a_0+jF_{k_*})\right\} = 
\left\{
\frac{F_{n-1}}{F_n}a_0 \right\} + O_\sigma \left(\frac{F_n^\sigma}{F_{k_*}^2}\right),
$$
and at the same time
$$
uF_n\le 
a_0+jF_{k_*}\le 
 (u+\eta/2)F_n
+
O_\sigma \left({F_n^\sigma}{F_{k_*}}\right)
.
$$
We should take into account that $F_{k_*} = O\left({F_n}^{\frac{1}{2}+\sigma}\right)$,
and everything is proven.$\Box$

\vskip+0.3cm
{\bf Acknowledgements.}

\vskip+0.3cm
The author is grateful to Nikita Shulga  and Johannes Schleischitz for the useful discussion.

The 
 research is supported by Austrian Science Fund (FWF), Forschungsprojekt PAT1961524.

\vskip+5cm

  \vskip+2cm
  
  author: Nikolay Moshchevitin,
  
 Technische Universität Wien
 
Institut für  diskrete Mathematik und Geomertie

Wiedner Hauptsrtaße  8

A-1040 Wien

Austria

  e-mails: nikolai.moshchevitin@tuwien.ac.at, moshchevitin@gmail.com

\end{document}